\documentstyle[12pt]{article}

\def\C{C\!\!\!\!I}
\def\D{{\cal D}}
\def\DD{{\bf D}}
\def\E{{\cal E}}
\def\EE{{\bf E}}
\def\F{{\cal F}}
\def\G{{\cal G}}
\def\H{{\cal H}}
\def\I{{\cal I}}

\def\M{{\cal M}}
\def\MM{{\bf M}}
\def\N{{\cal N}}
\def\O{{\cal O}}
\def\P{{\cal P}}

\def\rh{{\cal R\!\!H}}

\def\VV{{\bf V}}
\def\Z{Z\!\!\!Z}
\def\id{\mathop{\rm id}\nolimits}
\def\gr{\mathop{\rm gr}\nolimits}
\def\ker{\mathop{\rm Ker}\nolimits}
\def\sp{\mathop{\rm sp}\nolimits}
\def\spec{\mathop{\rm Spec}\nolimits}
\let\ov\overline

\newtheorem{theorem}{Theorem}[section]
\newtheorem{proposition}[theorem]{Proposition}
\newtheorem{lemma}[theorem]{Lemma}
\newtheorem{definition}[theorem]{Definition}
\newtheorem{corollary}[theorem]{Corollary}

\def\rem{\refstepcounter{theorem}\paragraph{Remark \thetheorem}}
\def\rems{\refstepcounter{theorem}\paragraph{Remarks \thetheorem}}
\def\proof{\paragraph{Proof}}
\def\inter#1{\subsection*{\hbox{}\hfill{\normalsize\sl #1}\hfill\hbox{}}}
\makeatletter \def\l@section{\@dottedtocline{1}{0em}{1.2em}} \makeatother

\begin{document}

\title{Moduli of pre-$\D$-modules, perverse sheaves\\ and the
Riemann-Hilbert morphism -I}
\author{Nitin Nitsure\thanks{Tata Institute of Fundamental
Research, Bombay} \and Claude Sabbah\thanks{CNRS, URA D0169,
Ecole Polytechnique, Palaiseau}} \date{March 28, 1995}
\maketitle

\begin{abstract} We construct a moduli scheme for semistable
pre-$\D$-modules with prescribed singularities and numerical
data on a smooth projective variety. These pre-$\D$-modules are
to be viewed as regular holonomic $\D$-modules with `level
structure'. We also construct a moduli scheme for perverse
sheaves on the variety with prescribed singularities and other
numerical data, and represent the de Rham functor (which gives
the Riemann-Hilbert correspondence) by an analytic morphism
between the two moduli schemes.
\end{abstract}

\vfill
\tableofcontents
\vfill
\newpage

\section{Introduction} This paper is devoted to the moduli
problem for regular holonomic $\D$-modules and perverse sheaves
on a complex projective variety $X$. It treats the case where
the singular locus of the $\D$-module is a smooth divisor $S$
and the characteristic variety is contained in the union of the
zero section $T^*_XX$ of the cotangent bundle of $X$ and the
conormal bundle $N^*_{S,X}$ of $S$ in $X$ (also denoted
$T_S^*X$). The sequel (part II) will treat the general case of
arbitrary singularities.

A moduli space for $\O$-coherent $\D$-modules on a smooth
projective variety was constructed by Simpson [S]. These are
vector bundles with integrable connections, and they are the
simplest case of $\D$-modules. In this moduli construction, the
requirement of semistability is automatically fulfilled by all
the objects.

Next in order of complexity are the so called `regular
meromorphic connections'. These $\D$-modules can be generated by
vector bundles with connections which have logarithmic
singularities on divisors with normal crossing. These
$\D$-modules are not $\O$-coherent, but are torsion free as
$\O$-modules. A moduli scheme does not exist for these
$\D$-modules themselves (see section 1 of [N]), but it is
possible to define a notion of stability and construct a moduli
for vector bundles with logarithmic connections. This was done
in [N]. Though many logarithmic connections give rise to the
same meromorphic connection, the choice of a logarithmic
connection is infinitesimally rigid if its residual eigenvalues
do not differ by nonzero integers (see section 5 of [N]).

In the present paper and its sequel, we deal with general
regular holonomic $\D$-modules. Such modules are in general
neither $\O$-coherent, nor $\O$-torsion free or pure
dimensional. We define objects called pre-$\D$-modules, which
play the same role for regular holonomic $\D$-modules that
logarithmic connections played for regular meromorphic
connections. We define a notion of (semi-)stability, and
construct a moduli scheme for (semi-) stable pre-$\D$-modules
with prescribed singularity stratification and other numerical
data. We also construct a moduli scheme for perverse sheaves
with prescribed singularity stratification and other numerical
data on a nonsingular variety, and show that the Riemann-Hilbert
correspondence defines an analytic morphism between (an open set
of) the moduli of pre-$\D$-modules and the moduli of perverse
sheaves.

The contents of this paper are as follows. Let $X$ be a smooth
projective variety, and let $S$ be a smooth hypersurface on $X$.
In section 2, we define pre-$\D$-modules on $(X,S)$ which may be
regarded as $\O_X$-coherent descriptions of those regular
holonomic $\D_X$-modules whose characteristic variety is
contained in $T^*_XX\cup T^*_SX$. The pre-$\D$-modules form an
algebraic stack in the sense of Artin, which is a property that
does not hold for the corresponding $\D$-modules.

In section 3, we define a functor from the pre-$\D$-modules to
$\D$-modules (in fact we mainly use the presentation of
holonomic $\D$-modules given by Malgrange [Mal], that we call
Malgrange objects). This is a surjective functor, and though not
injective, it has an infinitesimal rigidity property (see
proposition \ref{prop4}) which generalizes the corresponding
result for meromorphic connections.

In section 4, we introduce a notion of (semi-)stability for
pre-$\D$-modules, and show that semistable pre-$\D$-modules with
fixed numerical data form a moduli scheme.

Next, we consider perverse sheaves on $X$ which are
constructible with respect to the stratification $(X-S)\cup S$.
These have finite descriptions through the work Verdier,
recalled in section 5.

We observe that these finite descriptions are objects which
naturally form an Artin algebraic stack. Moreover, we show in
section 6 that S-equivalence classes (Jordan-H\"older classes)
of finite descriptions with given numerical data form a coarse
moduli space which is an affine scheme. Here, no hypothesis
about stability is necessary.

In section 7, we consider the Riemann-Hilbert correspondence.
When a pre-$\D$-module has an underlying logarithmic connection
for which eigenvalues do not differ by nonzero integers,
we functorially associate to it a finite description, which is
the finite description of the perverse sheaf associated to the
corresponding $\D$-module by the Riemann-Hilbert correspondence
from regular holonomic $\D$-modules to perverse sheaves. We show
that this gives an analytic morphism of stacks from the analytic
open subset of the stack (or moduli) of pre-$\D$-modules on
$(X,S)$ where the `residual eigenvalues are good', to the stack
(or moduli) of finite descriptions on $(X,S)$.

In section 8, we show that the above morphism of analytic stacks
is in fact a spread (surjective local isomorphism) in the
analytic category.  We also show that it has removable
singularities in codimension 1, that is, is can be defined
outside codimension two on any parameter space which is smooth
in codimension 1.

\paragraph{Acknowledgement} The authors thank the exchange
programme in mathematics of the Indo-French Center for the
Promotion of Advanced Research, New Delhi, the Ecole
Polytechnique, Paris, and the Tata Institute of Fundamental
Research, Bombay, for supporting their collaboration. The first
author also thanks ICTP Trieste and the University of
Kaiserslautern for their hospitality while part of this work was
done.

\section{Pre-$\D$-modules}

Let $X$ be a nonsingular variety and let $S\subset X$ be a
smooth divisor (reduced). Let $\I_S\subset \O_X$ be the ideal
sheaf of $S$, and let $T_X[\log S]\subset T_X$ be the sheaf of
all tangent vector fields on $X$ which preserve $\I_S$. Let
$\D_X[\log S]\subset \D_X$ be the algebra of all partial
differential operators which preserve $I_S$; it is generated as
an $\O_X$ algebra by $T_X[\log S]$.

The $\I_S$-adic filtration on $\O_X$ gives rise to a
(decreasing) filtration of $\D_X$ as follows: for $k\in\Z$
define $V^k\D_X$ as the subsheaf of $\D_X$ whose local sections
consist of operators $P$ which satisfy $P\cdot \I_S^j\subset
\I_S^{k+j}$ for all $j$. By construction, one has $\D_X[\log
S]=V^0\D_X$ and every $V^k(\D_X)$ is a coherent $\D_X[\log
S]$-module.

Let $p:N_{S,X}\to S$ denote the normal bundle of $S$ in $X$. The
graded ring $\gr_V\D_X$ is naturally identified with
$p_*\D_{N_{S,X}}$. Its $V$-filtration (corresponding to the
inclusion of $S$ in $N_{S,X}$ as the zero section) is then
split.

There exists a canonical section $\theta$ of the quotient ring
$\D_X[\log S]/\I_S\D_X[\log S]=\gr^0_V\D_X$, which is locally
induced by $x\partial_x$, where $x$ is a local equation for $S$.
It is a central element in $\gr^0_V\D_X$. This ring contains
$\O_S$ as a subring and $\D_S$ as a quotient (one has
$\D_S=\gr^0_V\D_X/\theta\gr^0_V\D_X$). One can identify locally
on $S$ the ring $\gr^0_V\D_X$ with $\D_S[\theta ]$.

A coherent $\gr^0_V\D_X$-module on which $\theta$ acts by $0$ is
a coherent $\D_S$-module. The locally free rank one
$\O_S$-module $\N_{S,X}=\O_X(S)/\O_X$ is a $\gr^0_V\D_X$-module
on which $\theta$ acts by $-1$.

\begin{definition}\rm
A {\sl logarithmic module} on $(X,S)$ will mean a sheaf of
$\D_X[\log S]$-modules, which is coherent as an $\O_X$-module. A
{\sl logarithmic connection} on $(X,S)$ will mean a logarithmic
module which is coherent and torsion-free as an $\O_X$-module.
\end{definition}

\rem
It is known that when $S$ is nonsingular, any logarithmic
connection on $(X,S)$ is locally free as an $\O_X$-module.

\begin{definition}[Family of logarithmic modules]\rm
Let $f:Z\to T$ be a smooth morphism of schemes. Let $Y\subset Z$
be a divisor such that $Y\to T$ is smooth. Let $T_{Z/T}[\log
Y]\subset T_{X/Y}$ be the sheaf of germs of vertical vector
fields which preserve the ideal sheaf of $Y$ in $\O_Z$. This
generates the algebra $\D_{Z/T}[\log Y]$. A family of
logarithmic modules on $Z/T$ is a $\D_{Z/T}[\log Y]$-module
which is coherent as an $\O_Z$-module, and is flat over $\O_T$.
When $f:Z\to T$ is the projection $X\times T\to T$, and
$Y=S\times T$, we get a {\sl family of logarithmic modules on
$(X,S)$ parametrized by $T$}.
\end{definition}

\rem
The restriction to $S$ of a logarithmic module is acted on by
$\theta$: for a logarithmic connection, this is the action of
the residue of the connection, which is an $\O_S$-linear
morphism.

\rem
There is an equivalence (restriction to $S$) between logarithmic
modules supported on the reduced scheme $S$ and
$\gr^0_V\D_X$-modules which are $\O_S$-coherent, (hence locally
free $\O_S$-modules, since they are locally $\D_S$-modules). In
the following, we shall not make any difference between the
corresponding objects.

\bigskip We give two definitions of pre-$\D$-modules. The two
definitions are `equivalent' in the sense that they give not
only equivalent objects, but also equivalent families, or more
precisely, the two definitions give rise to isomorphic algebraic
stacks. To give a familier example of such an equivalence, this
is the way how vector bundles and locally free sheaves are
`equivalent'. Note also that mere equivalence of objects is not
enough to give equivalence of families --- for example, the
category of flat vector bundles is equivalent to the category of
$\pi_1$ representations, but an algebraic family of flat bundles
gives in general only a holomorphic (not algebraic) family of
$\pi_1$ representations.

In their first version, pre-$\D$-modules are objects that live
on $X$, and the functor from pre-$\D$-modules to $\D$-modules
has a direct description in their terms. The second version of
pre-$\D$-modules is more closely related to the Malgrange
description of $\D$-modules and the Verdier description of
perverse sheaves, and the Riemann-Hilbert morphism to the stack
of perverse sheaves has direct description in its terms.

\begin{definition}[Pre-$\D$-module of first kind on $(X,S)$]\rm
Let $X$ be a nonsingular variety, and $S\subset X$ a smooth
divisor. A pre-$\D$-module $\EE = (E,F,t,s)$ on $(X,S)$ consists
of the following data

(1) $E$ is a logarithmic connection on $(X,S)$.

(2) $F$ is a logarithmic module on $(X,S)$ supported on the
reduced scheme $S$ (hence a flat connection on $S$).

(3) $t:(E\vert S) \to F$ and $s:F \to (E\vert S)$ are $\D_X[\log
S]$ linear maps,

which satisfies the following conditions:

(4) On $E\vert S$, we have $st = R$ where $R\in End(E\vert S)$
is the residue of $E$.

(5) On $F$, we have $ts = \theta_F$ where $\theta_F:F\to F$ is
the $\D_X[\log S]$ linear endomorphism induced by any Eulerian
vector field $x\partial /\partial x$.
\end{definition}

If $(E,F,t,s)$ and $(E',F',t',s')$ are two pre-$\D$-modules, a
morphism between them consists of $\D_X[\log S]$ linear
morphisms $f_0:E\to E'$ and $f_1:F\to F'$ which commute with
$t,t'$ and with $s,s'$.

\rem
It follows from the definition of a pre-$\D$-module $(E,F,t,s)$
that $E$ and $F$ are locally free on $X$ and $S$ respectively,
and the vector bundle morphisms $R$, $s$ and $t$ all have
constant ranks on irreducible components of $S$.

\paragraph{Example} Let $E$ be a logarithmic connection on
$(X,S)$. We can associate functorially to $E$ the following
three pre-$\D$-modules. Take $F_1$ to be the restriction of $E$
to $S$ as an $\O$-module. Let $t_1 = R$ (the residue) and $s_1 =
1_F$. Then $\EE_1=(E,F_1,t_1,s_1)$ is a pre-$\D$-module, which
under the functor from pre-$\D$-modules to $\D$-modules defined
later will give rise to the meromorphic connection corresponding
to $E$. For another choice, take $F_2 = E\vert S$, $t_2=1_F$ and
$s_2=R$. This gives a pre-$\D$-module $\EE_2 = (E,F_2,t_2,s_2)$
which will give rise to a $\D$-module which has nonzero torsion
part when $R$ is not invertible. For the third choice (which is
in some precise sense the minimal choice), take $F_3$ to be the
image vector bundle of $R$. Let $t_3 =R:(E\vert S)\to F_3$, and
let $s_3:F_3\hookrightarrow (E\vert S)$. This gives a
pre-$\D$-module $\EE_3 = (E,F_3,t_3,s_3)$. We have functorial
morphisms $\EE_3\to \EE_2 \to \EE_1$ of pre-$\D$-modules.

\begin{definition}[Families of pre-$\D$-modules]\rm
Let $T$ be a complex scheme. A family $\EE_T$ of
pre-$\D$-modules on $(X,S)$ parametrized by the scheme $T$, a
morphism between two such families, and pullback of a family
under a base change $T'\to T$ have obvious definitions (starting
from definition of families of $\D_X[\log S]$-modules), which we
leave to the reader. This gives us a fibered category $PD$ of
pre-$\D$-modules over the base category of $\C$ schemes. Let
$\cal PD$ be the largest (nonfull) subcategory of $PD$ in which
all morphisms are isomorphisms. This is a groupoid over $\C$
schemes.
\end{definition}

\begin{proposition}
The groupoid $\cal PD$ is an algebraic stack in the sense of
Artin.
\end{proposition}

\proof It can be directly checked that $\cal PD$ is a sheaf,
that is, descent and effective descent are valid for faithfully
flat morphisms of parameter schemes of families of
pre-$\D$-modules. Let $Bun_X$ be the stack of vector bundles on
$X$, and let $Bun_S$ be the stack of vector bundles on $S$. Then
$\cal PD$ has a forgetful morphism to the product stack
$Bun_X\times_{\C} Bun_S$. The later stack is algebraic and the
forgetful morphism is representable, hence the desired
conclusion follows.

\bigskip
Before giving the definition of pre-$\D$-modules of the second
kind, we observe the following.

\rem\label{rem1} Let $N$ be any line bundle on a smooth variety
$S$, and let $\ov{N} = P(N\oplus \O_S)$ be its projective
completion, with projection $\pi : \ov{N} \to S$. Let
$S^{\infty} = P(N)$ be the divisor at infinity. For any
logarithmic connection $E$ on $(\ov{N} ,S\cup S^{\infty})$, the
restriction $E\vert S$ is of course a $\D_{\ov{N}}[\log S\cup
S^{\infty}]$-module. But conversely, for any $\O $-coherent
$\D_{\ov{N}}[\log S\cup S^{\infty}]$-module $F$ scheme
theoretically supported on $S$, there is a natural structure of
a logarithmic connection on $(\ov{N} ,S\cup S^{\infty})$ on its
pullup $\pi ^*(F)$ to $\ov{N}$. The above correspondence is well
behaved in families, giving an isomorphism between the algebraic
stack of $\D_{\ov{N}}[\log S\cup S^{\infty}]$-modules $F$
supported on $S$ and the algebraic stack of logarithmic
connections $E$ on $(\ov{N} ,S\cup S^{\infty})$ such that the
vector bundle $E$ is trivial on the fibers of $\pi :\ov{N} \to
S$. The functors $\pi ^*(-)$ and $(-)\vert S$ are fully
faithful.

\rem\label{rem2} Let $S\subset X$ be a smooth divisor, and let
$N=N_{S,X}$ be its normal bundle. Then the following are
equivalent in the sense that we have fully faithful functors
between the corresponding categories, which give naturally
isomorphic stacks.

(1) $\D_X[\log S]$-modules which are scheme theoretically
supported on $S$.

(2) $\D_N[\log S]$-modules which are scheme theoretically
supported on $S$.

(3) $\D_{\ov{N}}[\log S \cup S^{\infty}]$-modules which are
scheme theoretically supported on $S$.

The equivalence between (2) and (3) is obvious, while the
equivalence between (1) and (2) is obtained as follows. The
Poincar\'e residue map $res:\Omega ^1_X[\log S] \to \O_S$ gives
the following short exact sequence of $\O_S$-modules. $$0\to
\Omega ^1_S \to \Omega ^1_X[\log S]|S \to \O_S\to 0$$ By taking
duals, this gives $$0 \to \O_S \to T_X[\log S]|S \to T_S\to 0.$$
It can be shown that there exists a unique isomorphism $T_X[\log
S]\vert S \to T_N[\log S]\vert S$ which makes the following
diagram commute, where the rows are exact.

$$\matrix{ 0 & \to & \O_S & \to & T_N[\log S]|S & \to & T_S &
\to & 0 \cr  & & \Vert & & \downarrow & & \Vert & & \cr 0 & \to
& \O_S & \to & T_X[\log S]|S & \to & T_S & \to & 0 \cr }$$

\rems
(1) Observe that the element $\theta$ is just the image of $1$
under the map $\O_S \to T_X[\log S]\vert S$.

(2) Using the notations of the  beginning of this section, one
can identify the ring $\pi_*\D_{\ov{N}}[\log S \cup S^{\infty}]$
with $\gr^0_V\D_X$. Hence $\theta$ is a global section of
$\D_{\ov{N}}[\log S \cup S^{\infty}]$.

\bigskip
We now make the following important definition.

\begin{definition}[Specialization of a logarithmic module]\rm
Let $E$ be a logarithmic module on $(X,S)$. Then the
specialization $\sp_SE$ will mean the logarithmic connection
$\pi ^*(E\vert S)$ on $(\ov{N_{S,X}} , S\cup S^{\infty})$.
\end{definition}

Now we are ready to define the second version of
pre-$\D$-modules.

\begin{definition}[Pre-$\D$-modules of the second kind on
$(X,S)$]\label{def1}\rm A pre-$\D$-mo\-dule (of the second kind)
$\EE = (E_0,E_1,c,v)$ on $(X,S)$ consists of the following data

(1) $E_0$ is a logarithmic connection on $(X,S)$,

(2) $E_1$ is a logarithmic connection on $(\ov{N_{S,X}},S\cup
S^\infty)$,

(3) $c:\sp_SE_0 \to E_1$ and $v:E_1 \to \sp_SE_0$ are
$\D_{\ov{N_{S,X}}}[\log S\cup S^\infty]$-linear maps,

which satisfies the following conditions:

(4) on $\sp_SE_0$, we have $cv = \theta_{\sp_SE_0}$,

(5) on $E_1$, we have $vc = \theta_{E_1}$,

(6) the vector bundle underlying $E_1$ is {\sl trivial} in the
fibers of $\pi:\ov{N_{S,X}}\to S$ (that is, $E_1$ is locally
over $S$ isomorphic to $\pi^*(E_1|S)$).
\end{definition}

If $(E_0,E_1,c,v)$ and $(E'_0,E'_1,c',v')$ are two
pre-$\D$-modules, a morphism between them consists of $\D_X[\log
S]$ linear morphisms $f_0:E_0\to E'_0$ and $f_1:E_1\to E'_1$
such that $\sp_Sf_0$ and $f_1$ commute with $v,v'$ and with
$c,c'$.

\begin{definition}[Families of pre-$\D$-modules of the second kind]\rm
Let $T$ be a complex scheme. A family $\EE_T$ of
pre-$\D$-modules on $(X,S)$ parametrized by the scheme $T$, a
morphism between two such families, and pullback of a family
under a base change $T'\to T$ have obvious definitions which we
leave to the reader. This gives us a fibered category $PM$ of
pre-$\D$-modules of second kind over the base category of $\C$
schemes.
\end{definition}

\begin{proposition}
The functor which associates to each family of pre-$\D$-module
$(E_0,E_1,c,v)$ of the second kind parametrized by $T$ the
family of pre-$\D$-module of the first kind $(E_0,E_1|S, c|S,
v|S)$ is an equivalence of fibered categories.
\end{proposition}

\proof This follows from remarks \ref{rem1} and \ref{rem2} above.

\section{From pre-$\D$-modules to $\D$-modules}

In this section we first recall the description of regular
holonomic $\D$-modules due to Malgrange [Mal] and we associate a
`Malgrange object' to a pre-$\D$-module of the second kind
(Proposition \ref{prop2}), which has good residual eigenvalues
(definition \ref{goodres}), each component of $S$ do not differ
by positive integers. Having such a direct description of the
Malgrange object enables us to prove that every regular
holonomic $\D$-module with characteristic variety contained in
$T^*_XX\cup T^*_SX$ arises from a pre-$\D$-module (Corollary
\ref{cor3}), and also helps us to prove an infinitesimal
rigidity property for the pre-$\D$-modules over a given
$\D$-module (Proposition
\ref{prop4}).

\inter{Malgrange objects}

Regular holonomic $\D$-modules on $X$ whose characteristic
variety is contained in $T^*_XX\cup T^*_SX$ have an equivalent
presentation due to Malgrange and Verdier, which we now
describe.

Let us recall the definition of the {\sl specialization}
$\sp_S(M)$ of a regular holonomic $\D_X$-module $M$: fix a
section $\sigma$ of the projection $\C\to\C/\Z$ and denote $A$
its image; every such module admits a unique (decreasing)
filtration $V^kM$ ($k\in\Z$) by $\D_X[\log S]$-submodules which
is good with respect to $V\D_X$ and satisfies the following
property: on $\gr^k_VM$, the action of $\theta$ admits a minimal
polynomial all of whose roots are in $A+k$. Then by definition
one puts $\sp_SM=\oplus_{k\in\Z}\gr^k_VM$. One has
$(\sp_SM)[*S]=\sp_S(M[*S])=(\gr_{V}^{\geq k}M)[*S]$ for all
$k\geq 1$, if we put $\gr_{V}^{\geq k}M=\oplus_{\ell\geq
k}\gr^\ell_VM$. The $p_*\D_{N_SX}$-module $\sp_SM$ does not
depend on the choice of $\sigma$ (if one forgets its gradation).

If $\theta$ acts in a locally finite way on a $\gr^0_V\D_X$ or a
$p_*\D_{N_{S,X}}$-module, we denote $\Theta$ the action of
$\exp(-2i\pi\theta)$.

Given a regular holonomic $\D_X$-module, we can functorially
associate to it the following modules:

(1) $M[*S]=\O_X[*S]\otimes_{\O_X}M$ is the $S$-localized
$\D_X$-module; it is also regular holonomic;

(2) $\sp_S M$ is the specialized module; this is a regular
holonomic $p_*\D_{N_SX}$-module, which is also {\sl monodromic},
i.e. the action of $\theta$ on each local section is locally (on
S) finite.

The particular case that we shall use of the result proved in
[Mal] is then the following:

\begin{theorem}
There is an equivalence between the category of regular
holonomic $\D_X$-modules and the category which objects are
triples $(\M,\overline M,\alpha)$, where $\M$ is a $S$-localized
regular holonomic $\D_X$-module, $\overline M$ is a monodromic
regular holonomic $p_*\D_{N_SX}$-module and $\alpha$ is an
isomorphism (of localized $p_*\D_{N_SX}$-modules) between
$\sp_S\M[*S]$ and $\overline M[*S]$.
\end{theorem}

In fact, the result of [Mal] does mention neither holonomicity
nor regularity. Nevertheless, using standard facts of the
theory, one obtains the previous proposition. Regularity
includes here regularity at infinity, i.e. along $S^\infty$.
This statement can be simplified when restricted to the category
of regular holonomic $\D$-modules which characteristic variety
is contained in the union $T^*_XX\cup T^*_SX$.

\begin{definition}\rm
A {\sl Malgrange object} on $(X,S)$ is a tuple $(M_0,M_1,C,V)$ where

(1) $M_0$ is an $S$-localized regular holonomic $\D_X$-module
which is a regular meromorphic connection on $X$ with poles on
$S$,

(2) $M_1$ is a $S$-localized monodromic regular holonomic
$p_*\D_{N_SX}$-module which is a regular meromorphic connection
on $N_{S,X}$ (or $\ov{N_{S,X}}$) with poles on $S$ (or on $S\cup
S^\infty$),

(3) $C,V$ are morphisms (of $p_*\D_{N_{S,X}}$-modules) between
$\sp_SM_0$ and $M_1$ satisfying $VC=\Theta-\id$ on
$\sp_SM_0$ and $CV=\Theta-\id$ on $M_1$.
\end{definition}

The morphisms between two Malgrange objects are defined in an
obvious way, making them an abelian category.

The previous result can be translated in the following way, using [Ve]:

\begin{corollary} There is an equivalence between the category
of regular holonomic $\D$-modules which characteristic variety
is contained in $T^*_XX\cup T^*_SX$ and the category of
Malgrange objects on $(X,S)$.
\end{corollary}

\inter{From pre-$\D$-modules to Malgrange objects}

\begin{definition}\label{goodres}\rm
(1) We say that a logarithmic connection $F$ on $(X,S)$ has
{\sl good residual eigenvalues} if for each connected component $S_a$
of the divisor $S$, the residual eigenvalues $(\lambda _{a,k})$
of $F$ along $S_a$ do not include a pair $\lambda _{a,i},\lambda
_{a,j}$ such that $\lambda _{a,i}-\lambda _{a,j}$ is a nonzero
integer.

(2) We say that a pre-$\D$-module $\EE =(E_0,E_1,s,t)$ has
{\sl good residual eigenvalues} if the logarithmic connection
$E_0$ has good residual eigenvalues as defined above.
\end{definition}

We now functorially associate a Malgrange object $\MM=\eta
(\EE)=(M_0,M_1,C,V)$ to each pre-$\D$-module $\EE =
(E_0,E_1,c,v)$ on $(X,S)$ with $E_0$ having good residual
eigenvalues.

\rem\label{rem3} By definition of a pre-$\D$-module it follows
that the nonzero eigenvalues of $\theta_a$ on $E_0|S_a$ (the
residue along $S_a$) are the same as the nonzero eigenvalues of
$\theta_a$ on $E_{1,a}$.

\begin{proposition}{\bf(The Malgrange object associated to a
pre-$\D$-module with good residual
eigenvalues)}\quad\label{prop2} Let $\EE=(E_0,E_1,c,v)$ be a
pre-$\D$-module on $(X,S)$ of the second kind (definition
\ref{def1}), such that $E_0$ has good residual eigenvalues. Let
$\eta (\EE ) = (M_0,M_1,C,V)$ where

(1) $M_0=E_0[*S]$,

(2) $M_1=E_1[*S]$,

(3) $C = c\circ
\displaystyle{e_{}^{-2i\pi\theta_{E_0}}- 1\over\theta_{E_0}}$.

(4) $V=v$

Then $\eta (\EE)$ is a Malgrange object, and $\eta$ is
functorial in an obvious way.
\end{proposition}

\proof Because $E_0$ has good residual eigenvalues, one can use
the filtration $V^kE_0[*S]$ $=I_{S}^{k}E_0\subset E_0[*S]$ in
order to compute $\sp_SE_0[*S]$. It follows that the
specialization of $E_0[*S]$ when restricted to $N_{S,X}-S$ is
canonically isomorphic to the restriction of $\sp_SE_0=\pi
^*(E_0\vert S)$ to $N_{X,S}-S$.

\inter{Essential surjectivity}

\begin{proposition}\label{prop3}
Every Malgrange object $(M_0,M_1,C,V)$ on $(X,S)$ can be
obtained in this way from a pre-$\D$-module.
\end{proposition}

\proof This follows from [Ve]: one chooses Deligne lattices in
$M_0$ and $M_1$. One uses the fact that every $\D$-linear map
between holonomic $\D$-modules is compatible with the
$V$-filtration, so sends the specialized Deligne lattice of
$M_0$ to the one of $M_1$. Moreover, the map $v$ can be obtained
from $V$ because the only integral eigenvalue of $\theta$ on the
Deligne lattice is $0$, so
$\displaystyle{e_{}^{-2i\pi\theta}- 1\over\theta}$ is
invertible on it.

The previous two propositions give the following.

\begin{corollary}\label{cor3} The functor from pre-$\D$-modules
on $(X,S)$ to regular holonomic $\D$-modules on $X$ with
characteristic variety contained in $T^*_XX\cup T^*_SX$ is
essentially surjective.
\end{corollary}

\inter{Infinitesimal rigidity}

For a regular holonomic $\D$-module $\MM$ with characteristic
variety $T^*_XX\cup T^*_SX$, there exist several nonisomorphic
pre-$\D$-modules $\EE$ which give rise to the Malgrange object
associated to $\MM$. However, we have the following
infinitesimal rigidity result, which generalizes the
corresponding results in [N].

\begin{proposition}[Infinitesimal rigidity]\label{prop4} Let
$T=\spec\displaystyle{\C [\epsilon ]\over (\epsilon ^2)}$. Let $\EE_T$
be a family of pre-$\D$-modules on $(X,S)$ parametrized by $T$.
Let the associated family $\MM_T$ of $\D$-modules on $X$ be
constant (pulled back from $X$). Let $\EE$ (which is the
specialization at $\epsilon =0$) be of the form $\EE =
(E,F,s,t)$ where along any component of $S$, no two distinct
eigenvalues of the residue of the logarithmic connection $E$
differ by an integer. Then the family $\EE_T$ is also constant.
\end{proposition}

\proof By [N], the family $E_{T}$ is constant, as well as the
specialization $\sp_SE_{T}$. As a consequence, the residue
$\theta_{E_T}$ is constant. Let us now prove that the family
$F_T$ is constant.

Let $S_a$ be a component of $S$ along which the only possible
integral eigenvalue of $\theta_E$ is $0$. Then it is also the
only possible integral eigenvalue of $\theta_F$ along $S_a$
because the generalized eigenspaces of $\theta_E$ and $\theta_F$
corresponding to a nonzero eigenvalue are isomorphic by $s$ and
$t$ (see remark \ref{rem3}). We also deduce from [N] that $F_T$
is constant as a logarithmic module along this component.

Assume now that $0$ is not an eigenvalue of $\theta_E$ along
$S_a$ but is an eigenvalue of $\theta_F$ along this component.
Then $\theta_F$ may have two distinct integral eigenvalues, one
of which is $0$. Note that, in this case, $\theta_E$ is an
isomorphism (along $S_a$), as well as $\theta_{E_T}$ which is
obtained by pullback from $\theta_E$. It follows that on $S_a$
we have an isomorphism $F_T\simeq E_T|S_a\oplus
\ker\theta_{F_T}$. Consequently $\ker\theta_{F_T}$ is itself a
family. It is enough to show that this family is constant. But
the corresponding meromorphic connection on $N_{S,X}^{}-S$ is
constant, being the cokernel of the constant map $C_T:M_{0T}\to
M_{1T}$. We can then apply the result of [N] because the only
eigenvalue on $\ker\theta_F$ is $0$.

The maps $s_T$ and $t_T$ are constant if and only if for each
component $S_a$ of $S$ and for some point $x_a\in S_a$ their
restriction to $F_T|{x_a}\times T$ and $E_T|{x_a}\times T$ are
constant. This fact is a consequence of the following lemma.

\begin{lemma} Let $E$ and $F$ be finite dimensional complex
vector spaces, and let $\theta_E\in End (E)$ and $\theta_F\in
End (F)$ be given. Let $V\subset W=Hom(F,E) \times Hom(E,F)$ be
the closed subscheme consisting of $(s,t)$ with $st=\theta_E$
and $ts=\theta_F$. Let $\phi :W\to W$ be the holomorphic map
defined by $$\phi (s,t) = (s, t {e^{st} -1 \over st}).$$ Then the
differential $d\phi$ is injective on the Zariski tangent space
to $V$ at any closed point $(s,t)$.
\end{lemma}

\proof Let $(a, b)$ be a tangent vector to $V$ at $(s,t)$. Then
by definition of $V$, we must have $at+sb=0$ and $ta+bs=0$.
Using $at+sb=0$, we can see that $d\phi (a, b) = (a, bf(st))$
where $f$ is the entire function on $End(E_0)$ defined by the
power series $(e^x-1)/x$. Suppose $(a,bf(st))=0$. Then $a=0$ and
so the condition $ta+bs=0$ implies $bs=0$. As the constant term
of the power series $f$ is $1$ and as $bs=0$, we have
$bf(st)=b$. Hence $b=0$, and so $d\phi$ is injective.

\section{Semistability and moduli for pre-$\D$-modules.}

In order to construct a moduli scheme for pre-$\D$-modules, one needs a
notion of semistability. This can be defined in more than one way.
What we have chosen below is a particularly simple and canonical
definition of semistability. (In an earlier version of this
paper, we had employed a definition of semistability in terms of
parabolic structures, in which
we had to fix the ranks of $s:E_1\to E_0|S$ and $t:E_0|S \to
E_1$ and a set of parabolic weights.)

Let $S_a$ be the irreducible
components of the smooth divisor $S\subset X$. For a pre-$\D$-module
$\EE =(E_0,E_1,s,t)$, we denote by $E_a$ the restriction of
$E_1$ to $S_a$, and we denote by $s_a$ and $t_a$ the
restrictions of $s$ and $t$.

\inter{Definition of semistability}

We fix an ample line bundle on $X$, and denote the resulting
Hilbert polynomial of a coherent sheaf $F$ by $p(F,n)$.
For constructing a moduli, we fix the Hilbert
polynomials of $E_0$ and $E_a$, which we denote by $p_0(n)$ and
$p_a(n)$.
Recall (see [S]) that an $\O _X$-coherent $\D _X[\log S]$-module $F$
is by definition {\sl semistable} if it is pure dimensional, and
for each $\O _X$ coherent
$\D _X[\log S]$ submodule $F'$, we have the inequality
$p(F',n)/rank (F') \le p(F,n)/rank (F)  $
for large enough $n$. We call $p(F,n)/rank (F)$ the {\sl
normalized Hilbert polynomial} of $F$.

\begin{definition}\rm
We say that the pre-$\D$-module $\EE$ is {\sl semistable} if
the $\D_X[\log S]$-modules $E_0$ and $E_a$ are semistable.
\end{definition}

\rems
(1) It is easy to prove that the semistability of the $\D
_X[\log S]$-module $E_a$ is equivalent to the semistability of
the logarithmic connection $\pi ^*_a(E_a)$ on $P(N_{S_a,X}\oplus
1)$ with respect to a natural choice of polarization.

(2) When $X$ is a curve, a pre-$\D$-module $\EE$ is semistable if
and only if the logarithmic connection $E_0$ on $(X,S)$ is
semistable, for then $E_1$ is always semistable.

(3) Let the dimension of $X$ be more than one.
Then even when a pre-$\D$-module $\EE$ is a pre meromorphic connection
(equivalently, when $s:E_1 \to E_0\vert S$ is an isomorphism),
the definition of semistability of $\EE$ does not
reduce to the semistability of the underlying logarithmic
connection $E_0$ on $(X,S)$. This is to be expected because we
do not fix the rank of $s$ (or $t$) when we consider families of
pre-$\D$-modules. Also note that even in dimension one,
meromorphic connections are not a good subcategory of the
abelian category of all regular holonomic $\D$-modules with
characteristic variety contained in $T^*_XX\cup T^*_SX$, in the
sense that a submodule or a quotient module of a meromorphic
connection is not necessarily a meromorphic connection.

\inter{Boundedness and local universal family}
We let the index $i$ vary over $0$ and over the $a$.

\begin{proposition}[Boundedness] Semistable pre-$\D$-modules
with given Hilbert po\-lynomials $p_i$ form a bounded set, that
is, there exists a family of pre-$\D$-modules parametrized by a
noetherian scheme of finite type over $\C$ in which each
semistable pre-$\D$-module with given Hilbert polynomials occurs.
\end{proposition}

\proof This is obvious as each $E_i$ (where $i=0,a$) being
semistable with fixed Hilbert polynomial, is bounded.

Next, we construct a local universal family. By boundedness,
there exists a positive integer $N$ such that for all $n\ge N$,
the sheaves $E_0(N)$ and $E_1(N)$ are generated by global
sections and have vanishing higher cohomology. Let $\Lambda =
D_X[\log S]$. Let $\O _X =\Lambda_0 \subset \Lambda_1 \subset
\cdots \subset \Lambda$ be the increasing filtration of
$\Lambda$ by the order of the differential operators. Note that
each $\Lambda_k$ is an $\O_X$ bimodule, coherent on each side.
Let $r$ be a positive integer larger than the ranks of the
$E_i$. Let $Q_i$ be the quot scheme of quotients
$q_i:\Lambda_r\otimes \O_X (-N)^{p_i(N)}\to\!\!\!\!\to E_i$
where the right $\O_X$-module structure on $\Lambda_r$ is used
for making the tensor product. Note that $G_i=PGL(p_i(N))$ has a
natural action on $Q_i$. Simpson defines a locally closed
subscheme $C_i\subset Q_i$ which is invariant under $G_i$, and
a local universal family $E$ of $\Lambda$-modules parametrized
by $C_i$ with the property that for two morphisms $T\to C_i$,
the pull back families are isomorphic over an open cover
$T'\to T$ if and only if the two morphisms define $T'$ valued
points of $C_i$ which are in a common orbit of $G_i(T')$.

On the product $C_0\times C_a$, consider the linear schemes
$A_a$ and $B_a$ which respectively correspond to
$Hom_{\Lambda}(E_1,E_0)$ and $Hom_{\Lambda}(E_0,E_1)$ (see Lemma
2.7 in [N] for the existence and universal property of such linear
schemes). Let $F_a$ be the fibered product of $A_a$ and $B_a$
over $C_0\times C_a$. Let $H_a$ be the closed subscheme of $F_a$
where the tuples $(q_0,q_1,t,s)$ satisfy $st=\theta$ and
$ts=\theta$. Finally let $H$ be the fibered product of the
pullbacks of the $H_a$ to $C= C_0 \times \prod_a C_a$. Note that
$H$ parametrizes a tautological family of pre-$\D$-modules on
$(X,S)$ in which every semistable pre-$\D$-module with given
Hilbert polynomials occurs.

The group $$\G = G_0 \times \prod_a (G_a \times GL(1))$$ has a
natural action on $H$, with
$$(q_0,q_a,t_a,s_a)\cdot (g_0,g_a,\lambda_a) =
(q_0g_0,q_ag_a,(1/\lambda_a)t_a,\lambda_a s_a)$$
It is clear from the definitions of $H$ and this action that two
points of $H$ parametrise isomorphic pre-$\D$-modules if and only
if they lie in the same $G$ orbit.

The morphism $H\to C\times \prod _aC_a$ is an affine morphism
which is $\G$-equivariant, and by Simpson's construction of
moduli for $\Lambda$-modules, the action of $\G$ on $C\times
\prod _aC_a$ admits a good quotient in the sense of geometric
invariant theory. Hence a good quotient $H//\G$ exists by
Ramanathan's lemma (see Proposition 3.12 in [Ne]), which by
construction and universal properties of good quotients
is the coarse moduli scheme of semistable pre-$\D$-modules with
given Hilbert polynomials.

Note that under a good quotient in the sense of geometric
invariant theory, two different orbits can in some cases get mapped
to the same point (get identified in the quotient).
In the rest of this section, we determine what are the closed
points of the quotient $H//\G$.

\rem For simplicity of notation, we assume in the rest of this
section that $S$ has only one connected component. It will be
clear to the reader how to generalize the discussion when $S$
has more components.

\inter{Reduced modules}

Assuming for simplicity that $S$ has only one connected
component, so that $\G = \H \times GL(1)$ where
$H=G_0 \times G_1$, we can construct the quotient
$H//\G$ in two steps: first we go modulo
the factor $GL(1)$, and then take the quotient of $R=H//GL(1)$
by the remaining factor $\H$. The following lemma is obvious.

\begin{lemma}\label{lem4.5}
Let $T$ be a scheme of finite type over $k$, and let $V\to T$
and $W\to T$ be linear schemes over $T$. Let $V\times W$ be
their fibered product (direct sum) over $T$, and let $V\otimes
W$ be their tensor product. Let $\phi :V\times W\to V\otimes W$
be the tensor product morphism. Then its schematic image
$D\subset V\otimes W$ is a closed subscheme which (i) parametrizes all
decomposable tensors, and (ii) base changes correctly.
Let $GL(1,k)$ act on $V\times W$ by the formula
$\lambda \cdot (v,w) = (\lambda v, (1/\lambda )w)$. Then $\phi
:V\times W\to D$ is a good quotient for this action.
\end{lemma}

\proof
The statement is local on the base, so we can assume
that (i) the base $T$ is an affine scheme, and (ii) both the linear
schemes are closed linear subschemes of trivial vector bundles
on the base, that is, $V\subset A^m_T$ and $W\subset A^n_T$ are
subschemes defined respectively by homogeneous linear equations
$f_p(x_i)=0$ and $g_q(y_j)=0$ in the coordinates $x_i$ on
$A^m_T$ and $y_j$ on $A^n_T$.
Let $z_{i,j}$ be the coordinates on $A^{mn}_T$, so that
the map $\otimes :A^m_T\times _T A^n_T \to A^{mn}$ sends
$(x_i,y_j) \mapsto (z_{i,j})$ where $z_{i,j}=x_iy_j$.
Its schematic image is the subscheme $C$ of $A^{mn}_T$
defined by the equations $z_{a,b}z_{c,d} - z_{a,d}z_{b,c} = 0$,
that is, the matrix $(z_{i,j})$ should have rank $ < 2$.
Take $D$ to be the subscheme of $C$ defined by the equations
$f_p(z_{1,j},\ldots ,z_{m,j}) = 0$ and
$g_q(z_{i,1},\ldots ,z_{i,n}) = 0$. Now the lemma \ref{lem4.5} follows
trivially from this local coordinate description.

\paragraph{}
The above lemma implies the following.
To get the quotient $H//GL(1)$, we just have to
replace the fibered product $A\times B$ over
$C_0\times C_1$ by the closed subscheme $Z\subset D\subset
A\otimes B$, where $D$ is the closed subscheme consisting
of decomposable tensors $u$, and $Z$ is the closed subscheme of
$D$ defined as follows. Let $\mu _0$ and $\mu _1$ be the
canonical morphisms from $A\otimes B$ to the linear schemes
representing $End_{\Lambda} (E_0|S)$ and $End_{\Lambda} (E_1)$
respectively. Then
$Z$ is defined to consist of all $u$ such that $\mu
_0(u)=\theta \in End _{\Lambda}(E_0|S) $ and $\mu _1(u) = \theta
\in End_{\Lambda}(E_1)$. There is a canonical
$GL(1)$ quotient morphism $A\times B \to D$ over
$C_0\times C_1$, which sends $(s,t)\mapsto u=s\otimes t$.
These give the $GL(1)$ quotient map $H\to Z$.
Note that the map $H\to C_0\times C_1$ is $\G$ equivariant, and the action
of $GL(1)$ on $C_0\times C_1$ is trivial, so we get a $\H$-equivariant map
$Z\to C_0\times C_1$.

In order to describe the identifications brought about by the
above quotient, we make the following definition.

\begin{definition}\rm
A {\sl reduced module} is a tuple $(E_0,E_1,u)$ where $E_0$ and
$E_1$ are as in a pre-$\D$-module, and $u\in
Hom_{\Lambda}(E_1,E_0|S)\otimes Hom_{\Lambda}(E_0,E_1)$ is a
decomposable tensor,
such that the canonical maps $\mu _0:Hom_{\Lambda}(E_1,E_0|S)\otimes
Hom_{\Lambda}(E_0,E_1) \to End_{\Lambda}(E_0|S)$ and
$\mu _1: Hom_{\Lambda}(E_1,E_0|S)\otimes Hom_{\Lambda}(E_0,E_1)
\to End_{\Lambda}(E_1)$,
both map $u$ to the endomorphism $\theta$ of the appropriate
module. In other words,
there exist $s$ and $t$ such that $(E_0,E_1,s,t)$ is a
pre-$\D$-module, and $u=s\otimes t$.
We say that the reduced module $(E_0,E_1,s\otimes t)$ is the associated
reduced module of the pre-$\D$-module $(E_0,E_1,s,t)$.
Moreover, we say that a reduced module is semistable if it is
associated to a semistable pre-$\D$-module.
\end{definition}

\begin{lemma} Let $V$ and $W$ are two vector spaces, $v,v'\in V$ and
$w,w'\in W$, then

(1) If $v\otimes w=0$ then $v=0$ or $w=0$.

(2) If $v\otimes w=v'\otimes w'\ne 0$, then there exists a
scalar $\lambda \ne 0$ such that $v=\lambda v'$ and $w =
(1/\lambda ) w'$.
\end{lemma}

\rem The above lemma shows that if $\EE$ and $\EE '$ are two
non-isomorphic pre-$\D$-modules whose associated reduced modules
are isomorphic, then
we must have $s\otimes t =s'\otimes t'=0$. In particular,
$\theta$ will act by zero on $E_0|S$ and also on $E_1$ for such
pre-$\D$-modules as $st=0$ and $ts=0$.

\inter{S-equivalence and stability}

\begin{definition}\rm
By a {\sl filtration} of a logarithmic connection $E$ we shall
mean an increasing filtration $E_p$ indexed by $\Z$ by subvector
bundles which are logarithmic connections.
Similarly, a filtration on a $\D _X[\log S]$-module $F$
supported on $S$ will mean a filtration of the vector bundle
$F\vert S$ by subbundles $F_p$ which are $\D _X[\log
S]$-submodules. A filtration of a
pre-$\D$-module $(E_0,E_1,s,t)$ is an increasing filtration
$(E_i)_p$ of the logarithmic connection $E_i$ ($i=0,1$) such
that $s$ and $t$ are filtered morphisms with respect to the
specialized filtration of $E_0$ and the filtration of $E_1$.
A filtration of a reduced module
$(E_0,E_1,u)$, with $u=s\otimes t$ where we take $s=0$ and $t=0$
if $u=0$, is a filtration of the pre-$\D$-module
$(E_0,E_1,s,t)$. We shall always assume that this filtration is
exhaustive, that is, $(E_i)_p=0$ for $p\ll0$ and $(E_i)_p=E_i$
for $p\gg0$. A filtration is {\sl nontrivial} if some $(E_i)_p$
is a proper subbundle of $E_i$ for $i=0$ or $1$.

\end{definition}
For a filtered pre-$\D$-module (or reduced module), each step of
the filtration as well as the graded object are pre-$\D$-modules
(or reduced modules).

\rem\label{deform}
There is a natural family $(\EE_\tau)_{\tau\in A^1}^{}$ of
pre-$\D$-modules or reduced modules parametrized by the affine
line $A^1=\spec\C[\tau]$, which fibre at $\tau=0$ is the graded
object $\EE'$ and the fibre at $\tau_0\neq0$ is isomorphic to
the original pre-$\D$-module or reduced module $\EE$: put (for
$i=0,1$) $\E_i=\oplus_{p\in\Z}^{}(E_i)_p\tau^p\subset E_i\otimes
\C[\tau,\tau_{}^{-1}]$ and the relative $\D\log$-structure is
the natural one.

\begin{definition}\rm
A {\sl special filtration of a coherent $\O _X$-module} $E$ is a
filtration for which each $E_p$ has the same normalized Hilbert
polynomial as $E$. A {\sl special filtration of a reduced
module} $(E_0,E_1,u)$ is a filtration of this reduced module
which is special on $E_0$ and on $E_1$.
\end{definition}

The graded reduced module $\EE'$ associated with a special
filtration of a semistable reduced module $\EE$ is again
semistable.

\begin{definition}\rm
The equivalence relation on the set of isomorphism classes
of all semistable reduced modules generated by this relation
(by which $\EE '$ is related to $\EE$) will be called S-equivalence.
\end{definition}

\begin{definition}\label{defstable}\rm
We say that a semistable reduced module is {\sl stable} if it
does not admit any nontrivial special filtration.
\end{definition}

\rems
(1) Note in particular that if each $E_0$, $E_a$ is stable
as a $\Lambda$-module, then the
reduced module $\EE $ is
stable. Consequently we have the following. Though the
definition of stability depends on the ample line bundle $L$ on
$X$, irrespective of the choice of the ample bundle, for any
pre-$\D$-module such that the monodromy representation of
$E_0\vert (X-S)$ is irreducible, and the monodromy
representation of $\pi _a ^*E_a \vert (N_{S_a,X}-S_a)$ is
irreducible for each component $S_a$, the corresponding reduced
module is stable. The converse is not true
-- a pre-$\D$-module, whose reduced module is stable, need not
have irreducible monodromies. The example 2.4.1 in [N] gives a
logarithmic connection, whose associated pre-$\D$-module in
which $s$ is identity and $t$ is the residue, gives a stable
reduced module, but the monodromies are not irreducible.

(2) If $u=0$, the reduced module is stable if and only $E_0$ and
each $E_a$ is stable.

(3) When $X$ is a curve, a reduced module with $u\ne 0$ is stable if
and only if the logarithmic connection $E_0$ is stable. If
$u=0$, each $E_a$ must moreover have length at most one as an
$\O_X$-module. Hence over curves, there is a plentiful supply of
stable reduced modules.

\begin{lemma}\label{uisflat}
Let $(E_0,E_1,u)$ be a reduced module and let $(E_i)_p$ be
filtrations of $E_i$ ($i=0,1$). Then $s$ and $t$ are filtered
morphisms with respect to the specialized filtrations if and only
if there exists some point $P\in S$ such that the restrictions
of $s$ and $t$ to the fibre $E_{i,P}$ at $P$ are filtered with
respect to the restricted filtrations.
\end{lemma}

\proof
This follows from the fact that if a section $\sigma$ of a
vector bundle with a flat connection has a value $\sigma(P)$ in
the fibre at $P$ of a sub flat connection, then it is a section
of this subbundle: we apply this to $s$ (resp. $t$) as a section
of $Hom((E_0)_{p|S},(E_1)_{|S})$ (resp.
$Hom((E_1)_{p|S},(E_0)_{|S})$).

\inter{A criterion for stability}
Let $\EE=(E_0,E_1,u=s\otimes t)$ be a reduced module. Assume
that we are given filtrations $0=F_0(E_i)\subset
F_1(E_i)\subset\cdots\subset F_{\ell_i}(E_i)=E_i$ of $E_i$
($i=0,1$) by vector subbundles which are $\D_X[\log
S]$-submodules.

For $j=0,\ldots,\ell_i$ let  $k(j)$ be the smallest $k$ such
that $s(\sp_SF_j(E_0))\subset F_k(E_1)$ and let $J(s)$ be the
graph of the map $j\to k(j)$. A {\sl jump point} is a point
$(j,k(j))$ on this graph such that $k(j-1)<k(j)$. Consider also
the set $G_s$ made by points under the graph: $G_s=\{ (j,k)\mid
k\leq k(j)\}$. For $t$ there is an equivalent construction: we
have a map $k\to j(k)$ and a set $G_t$ on the left of the graph
$I(t)$: $G_t=\{ (j,k)\mid j\leq j(k)\}$.

\begin{definition}\rm
$u=s\otimes t$ is {\sl compatible} with the filtrations if the
two sets $G_s$ and $G_t$ intersect at most at (common) jump
points (where if $u=0$, take $s=0$ and $t=0$).
\end{definition}

\begin{proposition}\label{nonstable}
Let $\EE=(E_0,E_1,u)$ be a semistable reduced module. The
following conditions are equivalent:

(1) $\EE$ is not stable,

(2) there exists a nontrivial special filtration $F_j(E_i)$
($j=0,\ldots\ell_i$) of each $E_i$ where all inclusions are
proper and $u$ is compatible with these filtrations.
\end{proposition}

\proof
$(1)\Rightarrow(2)$: If $\EE$ is not stable, we can find two
nontrivial special filtrations $(E_0)_p$ and $(E_1)_q$ such that
$s$ and $t$ are filtered morphisms. Let $p_j$ ($j=1,\ldots
,\ell_0$) be the set of jumping indices for $(E_0)_p$ and $q_k$
($k=1,\ldots ,\ell_1$) for $(E_1)_q$. For each $j_0$ and $k_0$
we have $j(k(j_0))\leq j_0$ and $k(j(k_0))\leq k_0$. We define
$F_j(E_0)=(E_0)_{p_j}$ and $F_k(E_1)=(E_1)_{q_k}$. We get
nontrivial filtrations of $E_0$ and $E_1$ where all inclusions
are proper. Moreover there cannot exist two distinct points
$(j_0,k(j_0))$ and $(j(k_0),k_0)$ with $j_0\leq j(k_0)$ and
$k_0\leq k(j_0)$ otherwise we would have $j_0\leq j(k_0)\leq
j(k(j_0))\leq j_0$ and the same for $k_0$ so the two points
would be the equal. Consequently $u$ is compatible with these
filtrations.

$(2)\Rightarrow(1)$: We shall construct a special filtration
$((E_0)_p,(E_1)_q)$ of
the reduced module from the filtrations $F_j(E_i)$ of each $E_i$.
Choose a polygonal line with only positive slopes, going through
each jump point of $G_s$ and for which each jump point of $G_t$
is on or above this line (see figure \ref{fig1}).
\setlength{\unitlength}{.5truecm}
\begin{figure}[htb]
\begin{center}
\begin{picture}(10,8)(0,0)
\put(0,0){\line(1,0){10}}
\put(0,0){\line(0,1){8}}
\put(9.5,-.7){$j$}
\put(-.5,7.5){$k$}
\put(3,4){\circle*{.2}}
\put(6,6){\circle*{.2}}
\put(8,8){\circle*{.2}}
\put(3,4){\line(1,0){3}}
\put(6,6){\line(1,0){2}}
\put(8,8){\line(1,0){2}}
\put(3,0){\line(0,1){4}}
\put(6,4){\line(0,1){2}}
\put(8,6){\line(0,1){2}}
\put(1,1){\circle{.2}}
\put(2,3){\circle{.2}}
\put(3,4){\circle{.2}}
\put(5,5){\circle{.2}}
\put(7,7){\circle{.2}}
\put(0,1){\line(1,0){1}}
\put(1,3){\line(1,0){1}}
\put(2,4){\line(1,0){1}}
\put(3,5){\line(1,0){2}}
\put(5,7){\line(1,0){2}}
\put(1,1){\line(0,1){2}}
\put(2,3){\line(0,1){1}}
\put(3,4){\line(0,1){1}}
\put(5,5){\line(0,1){2}}
\put(7,7){\line(0,1){1}}
\put(6,3){$G_s$}
\put(2,6){$G_t$}
\multiput(.2,.2)(.2,.2){4}{\circle*{.1}}
\multiput(1,1)(.2,.3){10}{\circle*{.1}}
\multiput(3,4)(.4,.2){5}{\circle*{.1}}
\multiput(5,5)(.2,.2){15}{\circle*{.1}}
\end{picture}
\caption{\label{fig1}$\bullet=$ jump points of $s$, $\circ=$
jump points of $t$}
\end{center}
\end{figure}
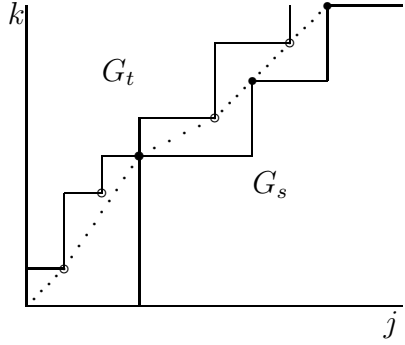
Choose increasing functions $p(j)$ and $q(k)$ such that
$p(j)-q(k)$ is identically $0$ on this polygonal line, is $<0$
above it and $>0$ below it (for instance, on each segment
$[(j_0,k_0),(j_1,k_1)]$ of this polygonal line, parametrised by
$j=j_0+m\varepsilon_1$, $k=k_0+m\varepsilon_2$, put
$p(j)=p(j_0)+\varepsilon_2(j-j_0)$ and
$q(k)=q(k_0)+\varepsilon_1(k-k_0)$, and $p(0)=q(0)=0$). For
$p(j)\leq p<p(j+1)$ put $(E_0)_p=F_j(E_0)$ and for $q(k)\leq
q<q(k+1)$ put $(E_1)_q=F_k(E_1)$. The filtration
$((E_0)_p,(E_1)_q,u)$ is then a nontrivial special filtration of
the reduced module $\EE$.

\begin{proposition}\label{open}
Semistability and stability are Zariski open conditions on the
parameter scheme of any family of reduced modules.
\end{proposition}

\proof As semistability is an open condition on $\D _X[\log
S]$-modules, it follows it is an open condition on reduced
modules. Now, for any family of semistable reduced modules
parametrised by a scheme $T$, all possible special filtrations
of the form given by \ref{nonstable}
on the specializations of the family are parametrised by
a scheme $U$ which is projective over $T$. The image of $U$ in
$T$ is the set of non stable points in $T$, hence its complement
is open.

\inter{Points of the moduli}

We are now ready to prove the following theorem.

\begin{theorem} Let $X$ be a projective variety together with an
ample line bundle, and let $S\subset X$ be a smooth divisor.

(1) There exists a coarse moduli scheme $\P$ for
semistable pre-$\D$-modules on $(X,S)$ with given Hilbert
polynomials $p_i$.  The scheme $\P$ is quasiprojective, in
particular, separated and of finite type over $\C$.

(2) The points of $\P$
correspond to S-equivalence classes of semistable
pre-$\D$-modules.

(3) The S-equivalence class of a semistable reduced module $\EE $
equals its isomorphism class if and only if $\EE $ is stable.

(4) $\P$ has an open subscheme $\P ^s$ whose points are the
isomorphism classes of all stable reduced modules. This is a
coarse moduli for (isomorphism classes of) stable reduced
modules.
\end{theorem}

\proof Let $\P = H//\G$. Then (1) follows by the construction of
$\P$. To prove (2), first note that by the existence of the
deformation $\EE _t$ (see \ref{deform}) of any reduced module
$\EE$ corresponding to a
weighted special filtration, and by the separatedness of $\P$,
the reduced module $\EE $ and its limit $\EE '$ go to the same point
of $\P$. Hence an S-equivalence class goes to a common point of
$\P$. For the converse, first recall that $\G = \H \times
GL(1)$, and the quotient $\P$ can
be constructed in two steps: $\P = R//\H$ where $R=H/\G$.
The scheme $R$ parametrizes a canonical family of reduced modules.
Let the $\H$ orbit of a point $x$ of $R$ corresponding the reduced
module $\EE $ not be closed in $R$. Let $x_0$ be any of its limit points.
Then there exists a 1-parameter subgroup $\lambda$ of $\H$ such
that $x_0 = \lim _{t\to 0} \lambda (t) x$. This defines a map
from the affine line $A^1$ to $R$, which sends $t\mapsto \lambda
(t)x$. Let $\EE _t$ be the pull back of the tautological family
of reduced modules parametrized by $R$.
Then from the description of the limits of the actions of
1-parameter subgroups on a quot scheme given in section 1 of
Simpson [S], it follows that $\EE$ has a special filtration such
that the family $\EE _t$ is isomorphic to a deformation of the type
constructed in \ref{deform} above. Hence the reduced modules
parametrized by $x$ and $x_0$ are S-equivalent. This proves (2).

If the orbit of $x$ is not closed, then it has a limit $x_0$
outside it under a 1-parameter subgroup, which by above represents a
reduced module $\EE '$ which is the limit of $\EE $ under a special
filtration. As by assumption $\EE '$ is not isomorphic to $\EE $, the
special filtration must be nontrivial. Hence $\EE $ is not stable.
Hence stable points have closed orbits in $R$. If $x$ represents
a stable reduced module, then $x$ cannot be the limit point of
any other orbit. For, if $x$ is a limit point of the orbit of
$y$, then by openness of stability (see \ref{open}), $y$ should
again represent a stable reduced module. But then by above, the
orbit of $y$ is closed. This proves (3).

Let $H^s\subset H$ be the open subscheme where the corresponding
pre-$\D$-module is stable. By (2) and (3) above, $H^s$ is
saturated under
the quotient map $H \to \P$, hence by properties of a good
quotient, its image $\P ^s$ is open in $\P$. Moreover by (2) and
(3) above, $H^s$ is the inverse image of $\P^ s$. Hence $H^s \to
\P ^s$ is a good quotient, which again by (2) and (3) is an
orbit space.  Hence points of $\P ^s$ are exactly the
isomorphism classes of stable reduced modules, which proves (4).

\section{Perverse sheaves, Verdier objects and finite descriptions}

Let $X$ be a nonsingular projective variety and let $S$ be a
smooth divisor. The abelian category of perverse sheaves
constructible with respect to the stratification $(X-S,S)$ of
$X$ is equivalent to the category of `Verdier objects' on
$(X,S)$. Before defining this category, let us recall the notion
of specialization along $S$.

Let $\E$ be a local system (of finite dimensional vector spaces)
on $X-S$. The {\sl specialization} $\sp_S\E$ is a local system
(of the same rank) on $N_{S,X}^{}-S$ equipped with an
endomorphism $\tau_\E$. It is constructed using the nearby cycle
functor $\psi$ defined by Deligne applied to the morphism which
describes the canonical deformation from $X$ to the normal
bundle $N_{S,X}^{}$.

A local system $\F$ on $N_{S,X}^{}-S$ equipped with an
endomorphism $\tau_\F$ is said to be {\sl monodromic} if
$\tau_\F$ is equal to the monodromy of $\F$ around $S$. Then
$\sp_S\E$ is monodromic.

\begin{definition}\rm
A {\sl Verdier object} on $(X,S)$ is a tuple $\VV=(\E,\F,C,V)$
where

(1) $\E$ is a local system on $X-S$,

(2) $\F$ is a monodromic local system on $N_{S,X}^{}-S$,

(3) $C:\sp_S\E\to\F$ and $V:\F\to\sp_S\E$ are morphisms of
(monodromic) local systems on $N_{S,X}^{}-S$ satisfying

(4) $CV=\tau_\F-\id$ and $VC=\tau_\E-\id$.
\end{definition}

\rem
The morphisms between Verdier objects on $(X,S)$ are defined in
an obvious way, and the category of Verdier objects is an
abelian category in which each object has finite length.
Hence the following definition makes sense.

\begin{definition}\rm
We say that two Verdier objects are {\sl S-equivalent} if they
admit Jordan-H\"older filtrations such that the corresponding
graded objects are isomorphic.
\end{definition}

\rem
Let $B$ be a tubular neighbourhood of $S$ in $X$, diffeomorphic
to a tubular neighbourhood of $S$ in $N_{S,X}^{}$. Put
$B^*=B-S$. The specialized local system $\sp_S\E$ can be
realized as the restriction of $\E$ to $B^*$, its monodromy
$\tau_\E$ at some point $x\in B^*$ being the monodromy along the
circle normal to $S$ going through $x$. Hence a Verdier object
can also be described as a tuple $\VV$ where $\F$ is a local
system on $B^*$ and $C$, $V$ are morphisms between $\E|B^*$ and
$\F$ subject to the same condition (4).

\bigskip
The notion of a family of perverse sheaves is not
straightforward. We can however define the notion of a family of
Verdier objects. Let us define first a family of local systems
on $X-S$ (or on $N_{S,X}^{}-S$) parametrized by a scheme $T$.
This is a locally free $p^{-1}\O_T$-module of finite rank, where
$p$ denotes the projection $X-S\times T\to T$. Morphisms between
such objects are $p^{-1}\O_T$-linear. The notion of a family of
Verdier objects is then straightforward.

In order make a moduli space for Verdier objects, we shall
introduce the category of `finite descriptions' on $(X,S)$. Let
us fix the following data (D):

(D1) finitely generated groups $G$ and $G_a$ for each component
$S_a$ of $S$,

(D2) for each $a$ an element $\tau_a$ which lies in the center
of $G_a$ and a group homomorphism $\phi_a:G_a\to G$.

\begin{definition}\label{def2}\rm
A finite description $\DD$ (with respect to the data (D)) is a
tuple $(E,\rho,F_a,\rho_a,C_a,V_a)$ where

(1) $\rho:G\to GL(E)$ is a finite dimensional complex
representation of  the group $G$; for each $a$ we will regard
$E$ as a representation of $G_a$ via the homomorphism
$\phi_a:G_a\to G$;

(2) for each $a$, $\rho_a:G_a\to GL(F_a)$ is a finite
dimensional complex representation of  the group $G$;

(3) for each $a$, $C_a:E\to F_a$ and $V_a: F_a\to E$ are
$G_a$-equivariant linear maps such that
$V_aC_a=\rho(\tau_a)-\id$ in $GL(E)$ and
$C_aV_a=\rho_a(\tau_a)-\id$ in $GL(F_a)$.
\end{definition}

A morphism between two finite descriptions has an obvious
definition.

\rem\label{rem4}
Let $P_0\in X-S$ and let $P_a$ be a point in the component
$B*_a$ of $B^*$. Choose paths $\sigma_a:[0,1]\to X-S$ with
$\sigma_a(0)=P_0$ and $\sigma_a(1)=P_a$. Let $G$ be the
fundamental group $\pi_1(X-S,P_0)$, and let $G_a = \pi_1
(B^*_a,P_a )$. Let $\tau_a \in G_a $ be the positive loop based
at $P_a $ in the fiber of $B^*_a\to S_a $. Finally, let
$\phi_a:G_a\to G$ be induced by the inclusion
$B^*_a\hookrightarrow X-S$ by using the path $\sigma_a$ to
change base points. Then, under the equivalence between
representations of fundamental group and local system, the
category of finite description with respect to the previous data
is equivalent to the category of Verdier objects on $(X,S)$.

\rem
The category of finite descriptions is an abelian category in
which each object has finite length. Therefore the notion of
S-equivalence as in definition 5.3 above makes sense for finite
descriptions.

\begin{definition}\rm
A family of finite descriptions parametrized by a scheme $T$ is
a tuple $(E_T, \rho_T,F_{T,a}, \rho_{T,a}, C_{T,a}, V_{T,a})$
where $E_T$ and the $F_{T,a}$ are locally free sheaves on $T$,
$\rho$ and $\rho_{T,a}$ are families of representations into
these, and the $C_{T,a}$ and $V_{T,a}$ are $\O_T$-homomorphisms
of sheaves satisfying the analogues of condition \ref{def2}.3
over $T$. The pullback of a family under a morphism $T'\to T$ is
defined in an obvious way, giving a fibered category. Let $PS$
denote the corresponding groupoid.
\end{definition}

\rem
It can be checked (we omit the details) that the groupoid $PS$
is an Artin algebraic stack.

\section{Moduli for perverse sheaves}

Let us fix data (D) as above.

\begin{theorem}
There exists an affine scheme of finite type over $\C$, which is
a coarse moduli scheme for finite descriptions
$\DD=(E,\rho,F_a,\rho_a,C_a,V_a)$ relative to {\rm (D)} with
fixed numerical data $n=\dim E$ and $n_a=\dim F_a$. The closed
points of this moduli scheme are the S-equivalence classes of
finite descriptions with given numerical data $(n,n_a)$.
\end{theorem}

Using remark \ref{rem4} we get

\begin{corollary}
There exists an affine scheme of finite type over $\C$, which is
a coarse moduli scheme for Verdier objects $\VV=(\E,\F,C,V)$ (or
perverse sheaves on $(X,S)$) with fixed numerical data $n={\rm
rank} \E$ and $n_a={\rm rank}\F|B^*_a$. The closed points of
this moduli scheme are the S-equivalence classes of Verdier
objects with given numerical data $(n,n_a)$.
\end{corollary}

The above corollary and its proof does not need $X$ to be a
complex projective variety, and the algebraic structure of $X$
does not matter. All that is needed is that the fundamental
group of $X-S$ and that of each $S_a$ is finitely generated.

The rest of this section contains the proof of the above
theorem.

\begin{proposition}\label{prop5}
(1) Let $\DD$ be a finite description, and let $\gr(\DD)$ be its
semisimplification. Then there exists a family $\DD_T$ of finite
descriptions parametrized by the affine line $T=A^1$ such that
the specialization $\DD_0$ at the origin $0\in T$ is isomorphic
to $\gr(\DD)$, while $\DD_t$ is isomorphic to $\DD$ at any $t\ne
0$.

(2) In any family of finite descriptions parametrized by a
scheme $T$, each S-equiva\-len\-ce class (Jordan-H\"older class)
is Zariski closed in $T$.
\end{proposition}

\proof
The statement (1) has a proof by standard arguments which we
omit. To prove (2), first note that if $\DD_T$ is any family and
$\DD'$ a simple finite description, then the condition that
$\DD' \times \{ t \}$ is a quotient of $\DD_t$ defines a closed
subscheme of $T$. From this, (2) follows easily.

\paragraph{Construction of Moduli} Let $E$ and $F_a$ be vector
spaces with $\dim(E)=n$ and $\dim(F_a) =n_a$. Let $\cal R$ be
the affine scheme of all representations $\rho$ of $G$ in $E$,
made as follows. Let $h_1,\ldots , h_r$ be generators of $G$.
Then $\cal R$ is the closed subscheme of the product $GL(E)^r$
defined by the relations between the generators. Similarly,
choose generators for each $G_a$, and let ${\cal R}_i$ be the
corresponding affine scheme of all representations $\rho_a$ of
$G_a$ in $F_a$.

Let $$A \subset {\cal R} \times \prod_a ({\cal R}_a \times
Hom(E,F_a) \times Hom(F_a,E))$$ be the closed subscheme defined
by condition \ref{def2}.3 above. Its closed points are tuples
$(\rho,\rho_a, C_a, V_a)$ where the linear maps $C_a:E\to F_a$
and $V_a:F_a\to E$ are $G_a$-equivariant under the
representations $\rho \phi_a: G_a\to GL(E)$ and $\rho_a: G_a\to
GL(F_a)$, and satisfy $V_aC_a = \rho (\tau_a )-1$ in $GL(E)$,
and $C_aV_a = \rho_a(\tau_a) -1$ in $GL(F_a)$ for each $a$.

The product group $\G =GL(E) \times (\prod_a GL(F_a))$ acts on
the affine scheme $A$ by the formula
$$(\rho ,\rho_a, C_a,V_a)\cdot (g,g_a) =
(g^{-1}\rho g, g_a^{-1}\rho_a g_a,g_a^{-1}C_ag, g^{-1}V_ag_a).$$
The orbits under this action are exactly the isomorphism classes
of finite descriptions.  The moduli of finite descriptions is
the good quotient $\F =A//\G$, which exists as $A$ is affine and
$\G$ is reductive. It is an affine scheme of finite type over
$\C$. It follows from \ref{prop5}.1 and \ref{prop5}.2 and
properties of a good quotient that the Zariski closures of two
orbits intersect if and only if the two finite descriptions are
S-equivalent. Hence closed points of $\F$ are S-equivalence
classes (Jordan-H\"older classes) of finite descriptions.

\section{Riemann-Hilbert morphism}

To any Malgrange object $\MM$, there is an obvious associated
Verdier object $\VV(\MM)$ obtained by applying the de~Rham
functor to each component of $\MM$. This defines a functor,
which is in fact an equivalence of
categories from Malgrange objects to Verdier objects. We have
already defined a functor $\eta$ from pre-$\D$-modules with good
residual eigenvalues to Malgrange objects. Composing, we get
an exact functor from pre-$\D$-modules with good residual eigenvalues
to Verdier objects. Choosing base points in $X$ and paths as in
remark \ref{rem4} we get an exact functor $\rh$ from pre-$\D$-modules
to finite descriptions.  This construction works equally well
for families of pre-$\D$-modules, giving us a holomorphic family
$\rh (\EE_T)$ of Verdier objects (or finite descriptions)
starting from a holomorphic family $\EE_T$ of pre-$\D$-modules
with good residual eigenvalues.

\rem
Even if $\EE_T$ is an algebraic family of pre-$\D$-modules with
good residual eigenvalues, the
associated family $\rh (\EE_T)$ of Verdier objects may not be
algebraic.

\rem If a semistable pre-$\D$-module has good residual
eigenvalues, then any other semistable pre-$\D$-module in its
S-equivalence class has (the same) good residual eigenvalues. Hence the
analytic open subset $T_g$ of the parameter space $T$ of any
analytic family of semistable pre-$\D$-modules defined by the
condition that residual eigenvalues are good is saturated under
S-equivalence.

\begin{lemma} If two semistable pre-$\D$-modules with good
residual eigenvalues are S-equivalent (in the sense of
definition \ref{defstable} above), then the associated finite
descriptions are S-equivalent (that is,
Jordan-H\"older equivalent).
\end{lemma}

\proof Let $\EE =(E_0,E_1,s,t)$ be a pre-$\D$-module with good
residual eigenvalues (that is, the logarithmic connection $E_0$
has good residual eigenvalues on each component of $S$) such that
$s\otimes t=0$. Then one can easily construct a family of
pre-$\D$-modules parametrized by the affine line $A^1$ which is the
constant family $\EE $ outside some point $P\in A^1$, and
specializes at $P$ to $\EE '=(E_0,E_1,0,0)$. Let $\phi :A^1 \to F$
be the resulting morphism to the moduli $\F$ of finite
descriptions. By construction, $\phi$ is constant on $A^1 -P$,
and so as $\F$ is separated, $\phi$ is constant. As the points
of $\F$ are the S-equivalence classes of finite descriptions, it
follows that the finite descriptions corresponding to $\EE $ and
$\EE '$ are S-equivalent.  Hence the S-equivalence class of the
finite description associated to a pre-$\D$-module depends only
on the reduced module made from the pre-$\D$-module. Now we must
show that any two S-equivalent (in the sense of \ref{defstable})
reduced semistable modules have associated finite descriptions
which are again S-equivalent (Jordan-H\"older equivalent). This
follows from the deformation given in \ref{deform} by using the
separatedness of $\F$ as above.

\bigskip
Now consider the moduli $\P = H//\G$ of semistable
pre-$\D$-modules. Let $H_g$ be the analytic open subspace of $H$
where the family parametrized by $H$ has good residual
eigenvalues. By the above remark, $H_g$ is saturated under $H\to
\P$. Hence its image $\P _g\subset \P$ is analytic open. Let
$\phi :H_g \to \F$ be the classifying map to the moduli $\F$ of
finite descriptions for the tautological family of
pre-$\D$-modules parametrized by $H$, which is defined because
of the the above lemma.  By the analytic universal property of
GIT quotients (see Proposition 5.5 of Simpson [S] and the remark
below), $\phi$ factors through an analytic map $\rh :P_g \to
\F$, which we call as the {\sl Riemann-Hilbert morphism}.

\rem In order to apply Proposition 5.5 of [S], note that
a $\G$-linear ample line bundle can be given on $H$ such that
all points of $H$ are semistable. Moreover, though the
proposition 5.5 in [S] is stated for semisimple groups, its proof
works for reductive groups.

\rem
The Riemann-Hilbert morphism can also be thought of as a
morphism from the analytic stack of pre-$\D$-modules with good
residual eigenvalues to the
analytic stack of perverse sheaves.

\section{Some properties of the Riemann-Hilbert morphism}

In this section we prove some basic properties of the morphism
$\rh$, which can be interpreted either at stack or at moduli
level.

\begin{lemma}[Relative Deligne construction]\label{lemdel}
(1) Let $T$ be the spectrum of an Artin local algebra of finite
type over $\C$, and let $\rho_T$ be a family of representations
of $G$ (the fundamental group of $X-S$ at base point $P_0$)
parametrized by $T$. Let $E$ be a logarithmic connection with
eigenvalue not differing by nonzero integers, such that the
monodromy of $E$ equals $\rho$, the specialization of $\rho_T$.
Then there exists a family $E_T$ of logarithmic connections
parametrized by $T$ such that $E_0=E$ and $E_T$ has monodromy
$\rho_T$.

(2) A similar statement is true for analytic germs of
$G$-representations.
\end{lemma}

\proof For each $a$, choose a fundamental domain $\Omega_a$ for
the exponential map ($z\mapsto \exp (2\pi \sqrt{-1} z)$) such
that the eigenvalues of the residue $R_a(E)$ of $E$ along $S_a$
are in the interior of the set $\Omega_a$. As the differential
of the exponential map $M(n,\C )\to GL(n,\C)$ is an isomorphism
at all those points of $M(n,\C )$ where the eigenvalues do not
differ by nonzero integers, using the fundamental domains
$\Omega_a$ we can carry out the Deligne construction locally to
define a family $E_T$ of logarithmic connections on $(X,S)$ with
$E_0=E$, which has the given family of monodromies.

Note that for the above to work, we needed the inverse function
theorem, which is valid for Artin local algebras.

\rem
If in the above, the family $\rho_T$ of monodromies is a
constant family (that is, pulled back under $T\to \spec (\C )$),
then $E_T$ is also a constant family as follows from Proposition
5.3 of [N].

\begin{proposition}[`Injectivity' of $\rh $]\label{propinj}
Let $\EE=(E,F,t,s )$ and $\EE'=(E',F',t',s')$ be
pre-$\D$-modules having good residual eigenvalues, such that for
each $a$, the eigenvalues of the residues of $E$ and $E'$ over
$S_a$ belong a common fundamental domain $\Omega_a$ for the
exponential map $exp :\C \to \C ^*:z\mapsto \exp (2\pi
\sqrt{-1}z)$. Then $\EE$ and $\EE'$ are isomorphic if and only
if the finite descriptions $\rh(\EE)$ and $\rh(\EE')$ are
isomorphic.
\end{proposition}

\proof
It is enough to prove that if the Malgrange objects $\MM$ and
$\MM'$ are isomorphic, then so are the pre-$\D$-modules $\EE$
and $\EE'$. First use the fact that, in a given meromorphic
connection $M$ on $X-S$ (or on $N_{S,X}^{}-S$), there exists one
and only one logarithmic connection having its residue along
$S_a$ in $\Omega_a$ for each $a$, to conclude that $E$ and $E'$
(resp. $F$ and $F'$) are isomorphic logarithmic modules. To
obtain the identification between $s$ and $s'$ (resp. $t$ and
$t'$), use the fact that these maps are determined by their
value at a point in each connected component $N_{S_a,X}^{}-S_a$
of $N_{S,X}^{}-S$ and this value is determined by the
corresponding $C_a$ or $C'_a$ (resp. $V_a$ or $V'_a$).

\begin{proposition}[Surjectivity of $\rh$]\label{propsurj}
Let $\DD$ be a finite description, and let $\sigma_a:\C ^*\to \C
$ be set theoretic sections of $z\mapsto \exp (2\pi
\sqrt{-1}z)$. Then there exists a pre-$\D$-module $\EE$ whose
eigenvalues of residue over $S_a$ are in image$(\sigma_a)$, for
which $\rh(\EE)$ is isomorphic to $\DD$.
\end{proposition}

\proof This follows from proposition \ref{prop3}.

\rem
The propositions \ref{propinj} and \ref{propsurj} together say
that the set theoretic fiber of $\rh$ over a given finite
description is in bijection with the choices of `good'
logarithms for the local monodromies of the finite description
(here `good' means eigenvalues do not differ by nonzero
integers).

\begin{proposition}[Tangent level injectivity for $\rh$]\label{propinfinj}
Let $(E,F,t,s)_T$ be a family of pre-$\D$-modules having good
residual eigenvalues parametrized by the spectrum $T$ of an
Artinian local algebra. Let the family $\rh(E,F,t,s)_T$ of
finite descriptions parametrized by $T$ be constant (pulled back
under $T\to \spec\C$). Then the family $(E,F,t,s)_T$ is
also constant.
\end{proposition}

\proof This is just the rigidity result of proposition \ref{prop4}.

\begin{proposition}[Infinitesimal surjectivity for $\rh$]\label{propinfsurj}
Let $T$ be the spectrum of an Artin local algebra of finite type
over $\C$, and let $\DD$ be a family of finite descriptions
parametrized by $T$. Let $\EE$ be a pre-$\D$-module having good
residual eigenvalues such that $\rh(\EE)=\DD_{\xi}$, the
restriction of $\DD$ over the closed point $\xi$ of $T$. Then
there exists a family $\EE'_T$ of pre-$\D$-modules having good
residual eigenvalues with $\EE'_{\xi}=\EE$ and
$\rh(\EE_T)=\DD_T$.
\end{proposition}

\proof This follows from lemma \ref{lemdel} and the proof of
proposition \ref{prop3} which works for families over Artin
local algebras.

\begin{theorem}
The analytic open substack of the stack (or analytic open subset
of the moduli) of pre-$\D$-modules on $(X,S)$, where $\EE $ has
good residual eigenvalues, is an analytic spread over the stack
(or moduli) of perverse sheaves on $(X,S)$ under the
Riemann-Hilbert morphism.
\end{theorem}

\proof This follows from propositions \ref{propsurj},
\ref{propinfinj} and \ref{propinfsurj} above.

Note that we have not defined $\rh$ on the closed analytic
subset $T_o$ of the parameter space of a family where $\EE $ does
not have good residual eigenvalues. Note that $T_o$ is defined
by a `codimension one' analytic condition, that is, if $T$ is
nonsingular, and if $T_o$ is a nonempty and proper subset of
$T$, then $T_o$ has codimension 1 in $T$. However, it follows
from Proposition \ref{propremov} below that the morphism $\rh$
on $T-T_o$ can be extended to an open subset of $T$ of
complementary codimension at least two. However, on the extra
points to which it gets extended, it may not represent the de
Rham functor.

\begin{proposition}[Removable singularities for $\rh$]\label{propremov}
Let $T$ be an open disk in $\C$ centered at $0$. Let
$\EE_T=(E,F,t,s)_T$ be a holomorphic family of pre-$\D$-modules
parametrized by $T$. Let the restriction $E_z$ have good
residual eigenvalues for all $z\in T-\{0\}$. Then there exists a
holomorphic family $\DD_U$ of finite descriptions parametrized
by a neighbourhood $U$ of $0\in T$ such that on $U-\{0\}$, the
families $\rh(\EE_U\vert U-\{0\}) $ and $\DD_{U- \{0\}}$ are
isomorphic.
\end{proposition}

If at $z=0$ the logarihmic connection $E$ does not have good
residual eigenvalues, it is possible to change it to obtain a
new logarithmic connection having good residual eigenvalues.
This is done by the classical `shearing transformation' that we
adapt below ({\sl inferior and superior modifications} for
pre-$\D$-modules). This can be done in family and has no effect
on the Malgrange object at least locally.

\begin{definition}\rm
If $E$ is a vector bundle on $X$, and $V$ a subbundle of the
restriction $E\vert S$, then the inferior modification ${_VE}$
is the sheaf of all sections of $E$ which lie in $V$ at points
of $S$. This is a locally free subsheaf of $E$ (but not
generally a subbundle). The superior modification $^VE$ is the
vector bundle $\O_X(S)\otimes {_VE}$.
\end{definition}

\rem\label{rem6}
If $E\vert S =V\oplus V'$, then we have a canonical isomorphism
$${_VE}\vert S \to V \oplus (\N^*_{S,X}\otimes V')$$ and hence
also a canonical isomorphism $$^VE|S\to (\N_{S,X}\otimes
V)\oplus V'$$

\rem\label{rem7}
If $(E,\nabla)$ is a logarithmic connection on $(X,S)$ and $V$
is invariant under the residue, then it can be seen that $_VE$
is invariant under $\nabla$, so is again a logarithmic
connection. We call it the inferior modification of the
logarithmic connection $E$ along the residue invariant subbundle
$V\subset E\vert S$. It has the effect that the residual
eigenvalues along $V$ get increased by $1$ when going from $E$
to $_VE$. As $\O_X(S)$ is canonically a logarithmic connection,
the superior modification $^VE$ is also a logarithmic
connection, with the residual eigenvalues along $V$ getting
decreased by $1$.

\bigskip
Let $(E,F,t,s)$ be pre-$\D$-module on $(X,S)$ such that $E$ has
good residual eigenvalues. Let us for simplicity of writing
assume that $S$ is connected. Let $E|S = \oplus_\alpha E^\alpha$
and $F=\oplus_\alpha F^\alpha$ be the respective direct sum
decompositions into generalized eigen subbundles for the action
of $\theta$. Then (see also remark \ref{rem3}) as $\theta$
commutes with $s$ and $t$, it follows that $t(E^\alpha) \subset
F^\alpha$ and $s(F_\alpha)\subset E^\alpha$. Moreover, when
$\alpha\ne 0$, the maps $s$ and $t$ are isomorphisms between
$E^\alpha$ and $F^\alpha$.

Now let $\alpha\ne 0$. Let $V=E^\alpha$ and $V'=\oplus_{\beta\ne
\alpha}E^\beta$. Let $F'' = \oplus_{\beta\ne \alpha}F^\beta$.
Let $F' = F^\alpha \oplus \N^*_{S,X}\otimes F''$. Let
$E'={_VE}$. Then using \ref{rem6} and the above, we get maps
$t':E'|S \to F'$ and $s':F'\to E'|S$ such that $(E',F',s',t')$
is a pre-$\D$-module.

\begin{definition}\label{definfmod}\rm
We call the pre-$\D$-module $(E',F',s',t')$ constructed above as
the inferior modification of $(E,F,s,t)$ along the generalized
eigenvalue $\alpha\ne 0$.
\end{definition}

Similarly, we can define the superior modification along a
generalized eigenvalue $\alpha\ne 0$ by tensoring with
$\O_X(S)$.

\rem
The construction of inferior or superior modification of
pre-$\D$-modules can be carried out over a parameter space $T$
(that is, for families) provided the subbundles $V$ and $V'$
form vector subbundles over the parameter space $T$ (their ranks
are constant).

\paragraph{Proof of \protect\ref{propremov}}
If the restriction $E= E_{T\vert z=0}$ has good residual
eigenvalues, then $\rh \EE_T$ has the desired property. So
suppose $E$ does not have good residual eigenvalues.

We first assume for simplicity of writing that $E$ fails to have
good residual eigenvalues because its residue $R_a$ on $S_a$ has
exactly one pair $(\alpha,\alpha-1)$ of distinct eigenvalues
which differ by a positive integer, with $\alpha-1\ne 0$. Let
$f_T$ be the characteristic polynomial of $R_{a,T}$. Then $f_0$
has a factorization $f_0 = gh$ such that the polynomials $g$ and
$h$ are coprime, $g(\alpha)=0$ and $h(\alpha-1)=0$. On a
neighbourhood $U$ of $0$ in $T$ we get a unique factorization
$f_T\vert U = g_Uh_U$ where $g_U$ specializes to $g$ and $h_U$
specializes to $h$ at $0$. By taking $U$ small enough, we may
assume that $g_U$ and $h_U$ have coprime specializations at all
points of $U$. Let $V_U$ be the kernel of the endomorphism
$g_U(R_{a,U})$ of the bundle $E_{a,U}$. If $U$ is small enough
then $F_U$ is a subbundle. Now take the inferior modification
$\EE'= ({_V}E_U, F'_U,t'_U,s'_U)$ of the family $(E,F,t,s)_U$ as
given by construction \ref{definfmod}. Then ${_VE}_U$ is a
family of logarithmic connections having good residual
eigenvalues, so by definition $\EE'$ has good residual
eigenvalues.

If $(0,1)$ are the eigenvalues, then use superior modification
along the eigenvalue $1$.

If $R_a$ has eigenvalues $(\alpha,\alpha-k)$ for some integer
$k\ge 1$, then repeat the above inferior (or superior)
modification $k$ times (whether to choose an inferior or
superior modification is governed by the following restriction :
the multiplicity of the generalized eigenvalue $0$ should not
decrease at any step). By construction, we arrive at the desired
family $(E',F',s',t')$.

\section*{References} \addcontentsline{toc}{section}{References}

[L] Laumon, G. : Champs alg\'ebriques. Preprint no. 88-33,
Universit\'e Paris Sud, 1988.

[Mal] Malgrange, B. : Extension of holonomic $\D$-modules, in
Algebraic Analysis (dedicated to M. Sato), M. Kashiwara and T.
Kawai eds., Academic Press, 1988.

[Ne] Newstead, P.E. : {\sl Introduction to moduli problems and
orbit spaces}, TIFR lecture notes, Bombay (1978).

[N] Nitsure, N. : Moduli of semistable logarithmic connections.
J. Amer. Math. Soc. 6 (1993) 597-609.

[S] Simpson, C. : Moduli of representations of the fundamental
group of a smooth projective variety - I, Publ. Math. I.H.E.S.
79 (1994) 47-129.

[Ve] Verdier, J.-L. : Prolongements de faisceaux pervers
monodromiques, Ast\'erisque 130 (1985) 218-236.

\bigskip
Addresses:

School of Mathematics, Tata Institute of Fundamental Research,
Homi Bhabha Road, Bombay 400 005, India. e-mail:
nitsure@tifrvax.tifr.res.in

Centre de Mathematiques, CNRS ura169, Ecole Polytechnique,
Palaiseau cedex, France. e-mail: sabbah@orphee.polytechnique.fr

\end{document}